\documentclass{article}
\usepackage{graphicx}
\usepackage{camnum}
\usepackage{harvard}

\renewcommand{\SS}{\CC{S}}

\newcommand{\ee}{{\mathrm e}}
\newcommand{\ii}{{\mathrm i}}
\newcommand{\hyper}[5]
       {{}_{#1} F_{#2} \!\left[
           \begin{array}{l}
         #3;\\#4;
           \end{array}#5\right]
       }

\title{On skyburst polynomials and their zeros}

\author{Mar\'{\i}a Jos\'e Cantero
\\Departamento de Matem\'atica Aplicada and IUMA\\ Escuela de Ingenier\'{\i}a y Arquitectura \\Universidad de Zaragoza\\ Spain\\ \&
\\Arieh Iserles \\
Department of Applied Mathematics and Theoretical Physics\\Centre for Mathematical Sciences\\University of Cambridge\\United Kingdom}

Â

\begin{document}
\maketitle
\tableofcontents

\begin{abstract}
 We consider polynomials orthogonal on the unit circle with respect to the complex-valued measure $z^{\omega-1}\D z$, where $\omega\in\BB{R}\setminus\{0\}$. We derive their explicit form, a generating function and several recurrence relations. These polynomials possess an intriguing pattern of zeros which, as $\omega$ varies, are reminiscent of a firework explosion. We prove this pattern in a rigorous manner.
\end{abstract}

\vspace{6pt}

\section{Introduction}

The subject matter of this paper is a specific family of polynomials, orthogonal on the unit circle with respect to the complex-valued measure $\D\mu(z)=z^{\omega-1}\D z$, where $\omega\in\BB{R}\setminus\{0\}$.

\begin{figure}[htb]
  \begin{center}
    \hspace*{-30pt}\includegraphics[width=130pt]{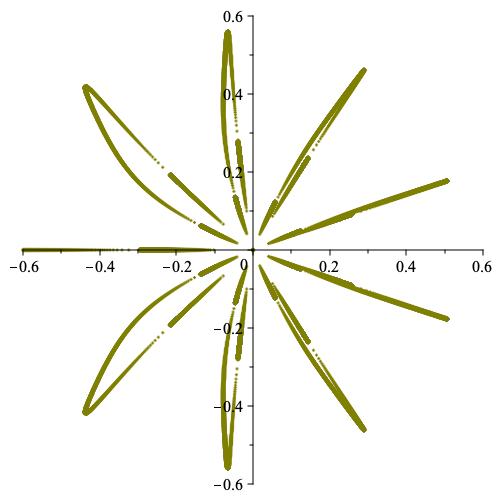}\hspace{10pt}\includegraphics[width=130pt]{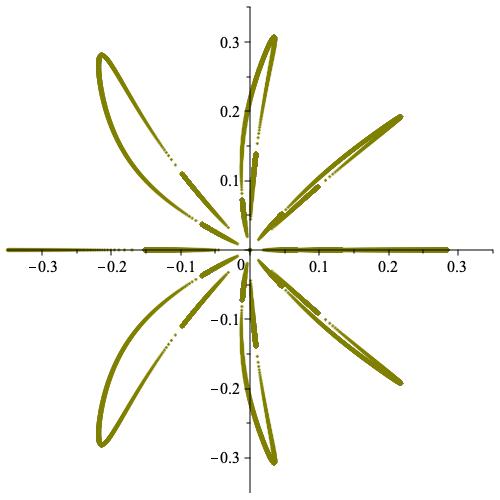}\hspace{10pt}\includegraphics[width=130pt]{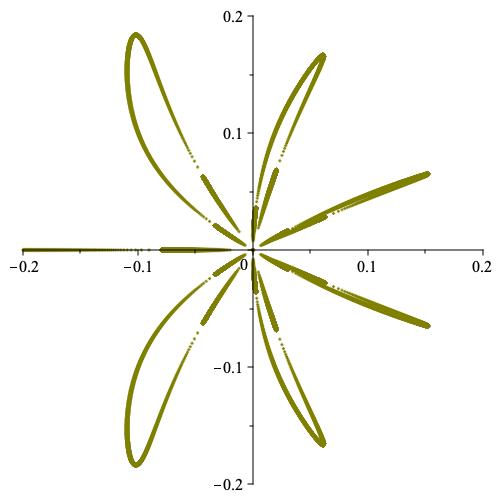}\hspace*{-30pt}
    
     \vspace*{12pt}
     \hspace*{-30pt}\includegraphics[width=130pt]{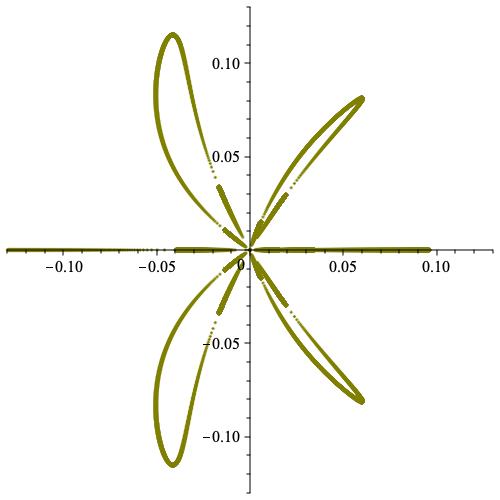}\hspace{10pt}\includegraphics[width=130pt]{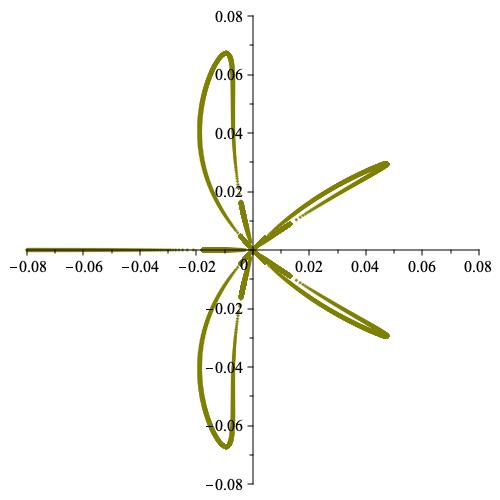}\hspace{10pt}\includegraphics[width=130pt]{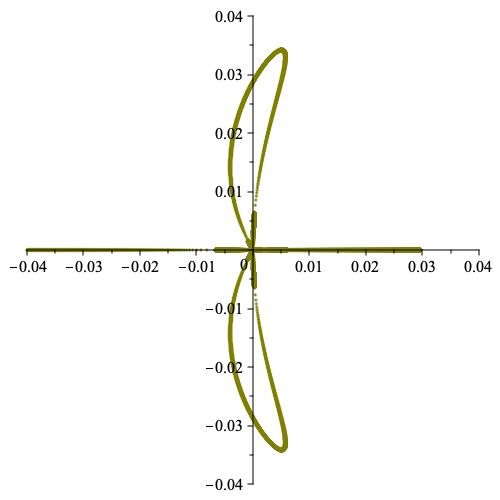}\hspace*{-30pt}
        
     \vspace*{12pt}
      \hspace*{-30pt}\includegraphics[width=130pt]{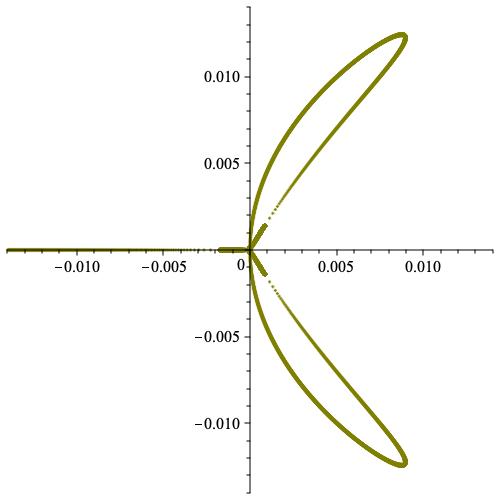}\hspace{10pt}\includegraphics[width=130pt]{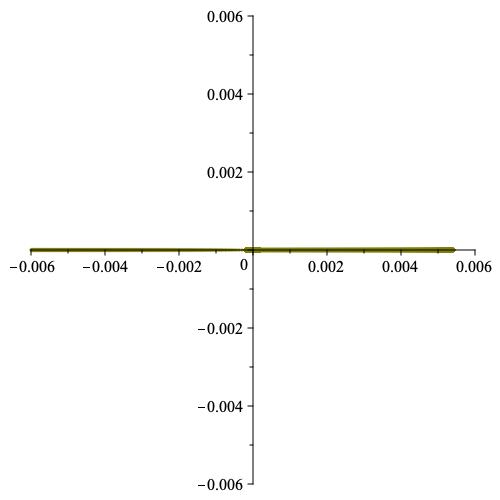}
    \caption{Zeros of $\CC{S}_9^\omega$ in intervals $m\leq\omega\leq m+1$, $m=0,1,\ldots,7$.}
    \label{Fig:1.1}
  \end{center}
  \begin{picture}(0,0)
     \put (26,462) {\small$\omega\in[0,1]$}
     \put (166,462) {\small$\omega\in[1,2]$}
     \put (307,462) {\small$\omega\in[2,3]$}
     \put (26,320) {\small$\omega\in[3,4]$}
     \put (166,320) {\small$\omega\in[4,5]$}
     \put (307,320) {\small$\omega\in[5,6]$}
     \put (82,177) {\small$\omega\in[6,7]$}
     \put (222,177) {\small$\omega\in[7,8]$}
  \end{picture}
\end{figure}
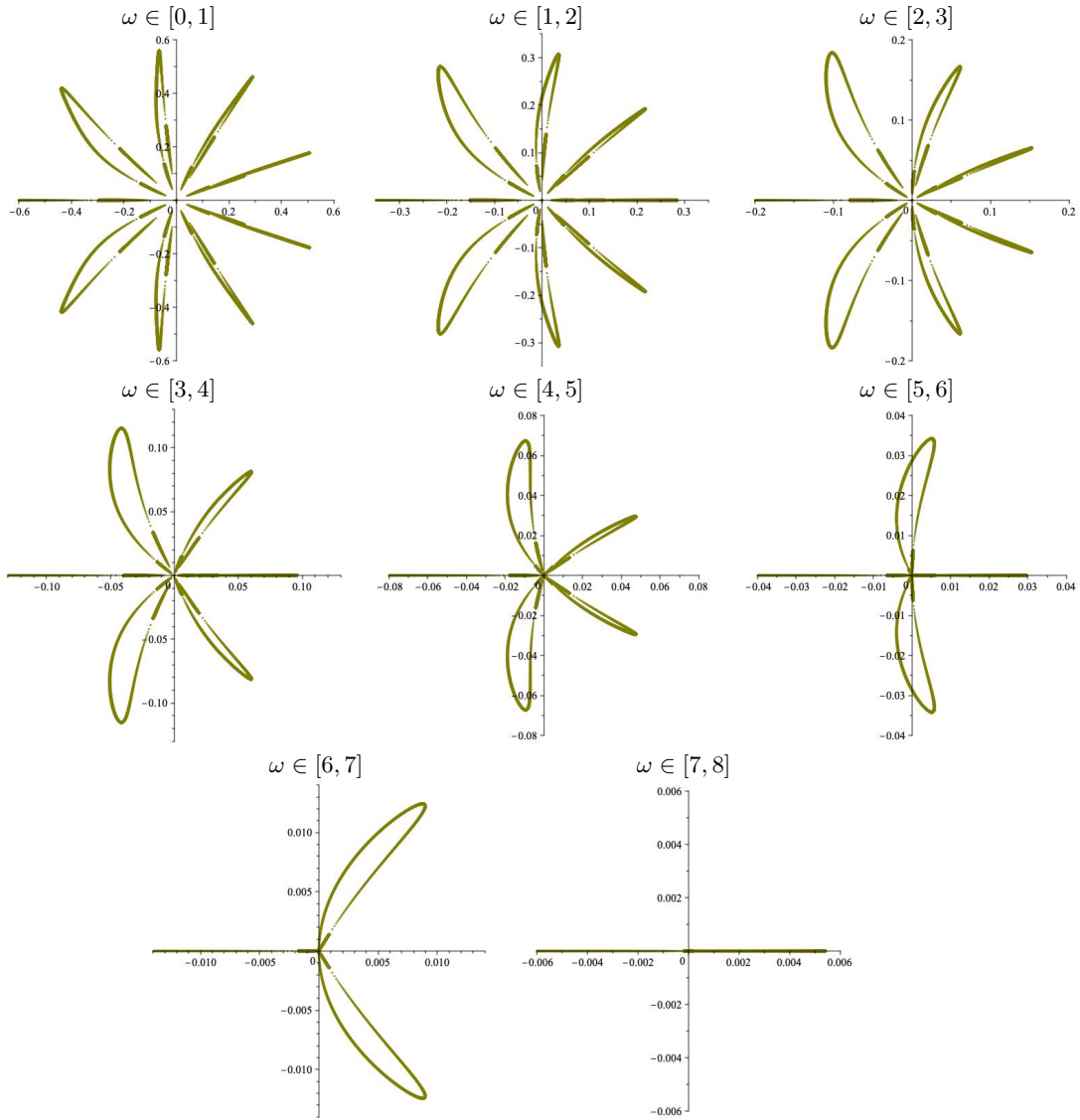

To motivate our interest (and the unusual name with which we have endowed them, {\em skyburst polynomials\/}), we commence in Fig.~\ref{Fig:1.1} with a plot of the zeros of $\SS_9^\omega$, the 9th-degree orthogonal polynomial with respect to the above measure, as $\omega\geq0$ increases. For $\omega=0$ the polynomial in question is $z^n$, with all its zeros at the origin. As $\omega$ grows in $(0,1)$, these zeros emerge from the origin into the complex plane: one into the interval $(-1,0)$, the remaining eight along equiangular rays in $\BB{C}$. Except for the zero in $(-1,0)$, as $\omega$ increases these zeros form eight loops, returning to the origin at $\omega=1$. Subsequently (in the second figure), the eight zeros emerge  for $\omega>1$ from the origin: one is `trapped' in $(-1,0)$, the rest form seven loops in the complex plane, eventually returning to the origin at $\omega=2$. Similar state of affairs unfolds in subsequent figures, each displaying zeros for $\omega\in[m,m+1]$ for increasing integer $m$: at each integer value of $m$ one more zero travels into $(-1,0)$ and stays there, the rest loop a loop in the complex plane, returning to the origin at the next integer value -- until $\omega=8$, when all the zeros live in $(-1,0)$ and remain there as $\omega>8$. 

Examining the trajectories of zeros as $\omega$ increases, the picture resembles a firework explosion, followed by a sequence of increasingly smaller (in both magnitude and complexity) explosions and eventually followed by a `fizzle'. This is the reason for the name ``skyburst polynomials''. 

A major feature of skyburst polynomials is that they possess a surprisingly simple explicit form: as proved in Section~2, the monic polynomials are 
\begin{displaymath}
  \SS_n^\omega(z)=z^n \hyper{2}{1}{-n,-\omega}{-n-\omega}{-\frac{1}{z}}.
\end{displaymath}
This, in itself, is fairly remarkable, because so few orthogonal polynomial systems on the unit circle are known in their explicit form \cite{cantero16fop,simon05opuc1} and it acts as a gateway towards a surprisingly simple generating function and a wealth of recurrence relations: this is the theme of Section~3. More effort is required to examine in great detail the observations that we have just made in Fig.~\ref{Fig:1.1} and prove them rigorously for every $n\in\BB{N}$. This is done in Section~4, building upon the material of Section~3.

In general, exceedingly little is known on orthogonal polynomials (whether on the real line or the unit circle) with respect to complex-valued measures. The one substantive result, \cite{celsus22kpt}, indicates that, for a specific, parameter-dependent family of polynomials orthogonal on the real line, the zeros behave in a highly complicated manner. That paper also emphasises the advantage of examining  the zeros as a function of a parameter: for each individual value we have just a number of points in $\BB{C}$ but, once we examine the evolution of these points as the parameter varies, the full (and intricate!) picture emerges. The current paper, the first to consider orthogonal polynomials on the unit circle in this setting, highlights a similar state of affairs. A snapshot of the zeros of $\SS_9^n$ (cf.~Fig.~\ref{Fig:1.1}) tells us very little but the complete `movie' as $\omega$ evolves indicates the underlying structure which, as we do in this paper, needs be subjected to rigorous analysis.

\setcounter{equation}{0}
\setcounter{figure}{0}
\section{Skyburst polynomials}

Let $\omega\in\BB{R}\setminus\{0\}$ and $\D\mu(z)=z^{\omega-1}\D z$. We consider the complex-valued bilinear form
\begin{equation}
  \label{eq:2.1}
  \langle f,g\rangle_\omega=\frac{1}{2\pi\ii}\int_{\bb{T}} f(z)\overline{g(\bar{z})} z^{\omega-1}\D z=\frac{1}{2\pi} \int_{-\pi}^\pi f(\ee^{\ii\theta}) \overline{g(\ee^{-\ii\theta})} \ee^{\ii\omega\theta}\D\theta,
\end{equation}
where $\BB{T}$ is the complex unit circle. Since $\omega\neq0$, \R{eq:2.1} is not a genuine inner product, since $\langle f,f\rangle_\omega=0$ does not imply $f=0$, yet this need not preclude the existence of orthogonal polynomials. Like in \cite{celsus22kpt} in the case of orthogonal polynomials on the real line with complex-valued measures, while important elements of classical theory are lost, the construct itself is amenable to analysis. While we cannot be assured {\em a priori\/}  for all $n\in\BB{Z}_+$ and $\omega$ of the existence of  a monic $n$th-degree polynomial $\SS_n^\omega$, orthogonal with respect to \R{eq:2.1}, it turns out that it always exists except for a finite number of values of $\omega$ and, moreover, such polynomials can be described explicitly and their zeros display intriguing patterns as $\omega$ varies.

We assume for the time being that $\omega>0$: this assumption, rendering our analysis considerably simpler, will be lifted toward the end of this section.

For reasons already explained in Section~1 and examined rigorously in the sequel, we call the $\SS_n^\omega$s {\em skyburst polynomials.\/}

We recall from \cite[p.65]{simon05opuc1} that an $n$th degree monic polynomial $p_n$, orthogonal on the unit circle with respect to the measure $\D\mu(z)=w(z)\D z$, can be represented in determinantal form,
\begin{equation}
  \label{eq:2.2}
  p_n(z)=\frac{1}{\GG{t}_n} \det\!\left[
  \begin{array}{ccccc}
    \mu_0 & \mu_1 & \cdots & \mu_{n-1} & 1\\
    \mu_{-1} & \mu_0 & \cdots & \mu_n & z\\
    \vdots & \vdots & & \vdots & \vdots\\
    \mu_{-n} & \mu_{-n+1} & \cdots & \mu_{-1} & z^n
  \end{array}
  \right]\!,
\end{equation}
using the moments of the underlying measure,
\begin{displaymath}
  \mu_n=\frac{1}{2\pi\ii}\int_{\bb{T}} z^n w(z)\frac{\D z}{z}=\frac{1}{2\pi} \int_{-\pi}^\pi \ee^{\ii n\theta} w(\ee^{\ii\theta})\D\theta,\qquad n\in\BB{Z}
\end{displaymath}
and
\begin{equation}
  \label{eq:2.3}
  \GG{t}_n=\left[
  \begin{array}{cccc}
  \mu_0 & \mu_1 & \cdots & \mu_{n-1}\\
  \mu_{-1} & \mu_0 & \cdots & \mu_{n-2}\\
  \vdots & \vdots & & \vdots\\
  \mu_{-n+1} & \mu_{-n+2} & \cdots & \mu_0
  \end{array}
  \right]\!,\qquad n\in\BB{Z}_+.
\end{equation}
Being purely algebraic constructs, \R{eq:2.2} and \R{eq:2.3} remain valid for complex-valued measures.

In the case of skyburst polynomials the moments are
\begin{displaymath}
  \mu_n^\omega=\frac{(-1)^n}{\pi} \frac{\sin\pi\omega}{n+\omega},\qquad n\in\BB{Z}, \quad \omega>0.
\end{displaymath}
Because of \R{eq:2.3}, $\SS_n^\omega$ exists if and only if $\GG{t}_n^\omega=\GG{t}_n\neq0$ and is bounded. 

\begin{lemma}
  Let $\omega\in\BB{R}\setminus\BB{Z}$, Then it is true that
  \begin{equation}
    \label{eq:2.4}
    \GG{t}_n^\omega= \left(\frac{\sin\pi\omega}{\pi\omega}\right)^{\!n} \frac{\prod_{\ell=0}^{n-1}\ell!^2}{\prod_{k=1}^{n-1}(k^2-\omega^2)^{n-k}},\qquad n\in\BB{Z}_+,\quad \omega\in\BB{R},
  \end{equation}
  hence $\SS_n^\omega$ exists.
\end{lemma}

\begin{proof}
  The moments of \R{eq:2.1} are
  \begin{displaymath}
    \mu_n^\omega=(-1)^n\frac{\sin\pi\omega}{\pi} \frac{1}{n+\omega},\qquad n\in\BB{Z},
  \end{displaymath}
  and substitution into \R{eq:2.3} results in
  \begin{Eqnarray*}
    \GG{t}_n&=&
    \left(\frac{\sin\pi\omega}{\pi}\right)^{\!n} \det\!\left[
    \begin{array}{ccccc}
       \frac{1}{\omega} & -\frac{1}{\omega+1} & \frac{1}{\omega+2} & \cdots & \frac{(-1)^{n-1}}{\omega+n-1}\\[2pt]
       -\frac{1}{\omega-1} & \frac{1}{\omega} & -\frac{1}{\omega+1} & \cdots & \frac{(-1)^{n-2}}{\omega+n-2}\\[2pt]
       \vdots & \vdots & \vdots & & \vdots\\[2pt]
       \frac{(-1)^{n-1}}{\omega-n+1} & \frac{(-1)^{n-2}}{\omega-n+2} & \frac{(-1)^{n-3}}{\omega-n+3} & \cdots & \frac{1}{\omega}
    \end{array}
    \right]\\
    &=&\left(\frac{\sin\pi\omega}{\pi}\right)^{\!n} \det\!\left[
  \begin{array}{ccccc}
  \frac{1}{\omega} & \frac{1}{\omega+1} & \frac{1}{\omega+2} & \cdots & \frac{1}{\omega+n-1}\\[2pt]
  \frac{1}{\omega-1} & \frac{1}{\omega} & \frac{1}{\omega+1} & \cdots & \frac{1}{\omega+n-2}\\[2pt]
  \vdots & \vdots & \vdots & & \vdots\\[2pt]
  \frac{1}{\omega-n+1} & \frac{1}{\omega-n+2} & \frac{1}{\omega-n+3} & \cdots & \frac{1}{\omega}
  \end{array}
  \right]\!.
  \end{Eqnarray*}
  The way we have obtained the second determinant is by pulling a factor of $-1$ from every odd (counting from zero) row and column of the first determinant. This gives a factor of $(-1)^{2n}=1$.
  
  We identify the last expression as a determinant of a {\em Cauchy matrix,\/} where the $(k,\ell)$ element is $1/(x_k+y_\ell)$, $k,\ell=0,\ldots,n-1$. The determinant of such a matrix is 
  \begin{displaymath}
    \frac{\prod_{k=1}^{n-1}\prod_{\ell=0}^{k-1}(x_k-x_\ell)(y_k-y_\ell)}{\prod_{k=0}^{n-1}\prod_{\ell=0}^{n-1}(x_k+y_\ell)}
  \end{displaymath}
  \cite{Schechter59oic}. The lemma follows by straightforward algebra once we let $x_k=\omega+k$, $y_k=-k$.
\end{proof}

\begin{corollary}
  The polynomial $\SS_n^\omega$ exists and is of degree $n$ for all $\omega\not\in\BB{N}$.
\end{corollary}

\begin{proof}
  Follows from \R{eq:2.4} because
  \begin{displaymath}
    \GG{t}_n^m=\frac{\prod_{\ell=0}^{n-1}\ell!^2}{(\pi m)^n} \prod_{\stackrel{\scriptstyle k=1}{k\neq m}}^{n-1} (k^2-m^2)^{k-n} \lim_{\omega\rightarrow m}\frac{\sin^n\pi\omega}{(m^2-\omega^2)^{n-m}}=0,\qquad m\in\BB{N},
  \end{displaymath}
  otherwise $\GG{t}_n^\omega\neq0$.
  
  Moreover, the coefficient of $z^n$ in $\SS_n^\omega$ is $\GG{t}_{n-1}/\GG{t}_n$, according to \R{eq:2.2}, and this is nonzero for $\omega\not\in\BB{N}$ by a similar argument.
\end{proof}

\begin{theorem}
  \label{th:explicit_form}
  Let $\omega>0$. It is true that
  \begin{equation}
    \label{eq:2.5}
    \SS_n^\omega(z)=z^n \hyper{2}{1}{-n,-\omega}{-n-\omega}{-\frac{1}{z}}\!,\qquad n\in\BB{Z}_+.
  \end{equation}  
\end{theorem}

\begin{proof}
  It is enough to prove that the monic polynomial given in \R{eq:2.5} is orthogonal to $z^k$, $k=0,\ldots,n-1$, with respect to the underlying bilinear form $\langle\,\cdot\,,\,\cdot\,\rangle_\omega$. It follows from \R{eq:2.5} that
  \begin{Eqnarray*}
     \langle\SS_n^\omega,z^k\rangle_\omega&=&\frac{\ii}{2\pi} \int_{-\pi}^\pi \sum_{\ell=0}^n (-1)^\ell \frac{(-n)_\ell(-\omega)_\ell}{\ell!(-n-\omega)_\ell} \ee^{\ii(n-\ell-k+\omega)\theta}\D\theta\\
     &=&\frac{1}{2\pi} \sum_{\ell=0}^n (-1)^\ell \frac{(-n)_\ell(-\omega)_\ell}{\ell!(-n-\omega)_\ell} \frac{\ee^{\ii(n-\ell-k+\omega)\pi}-\ee^{-\ii(n-\ell-k+\omega)\pi}}{n-\ell-k+\omega}\\
     &=&\frac{(-1)^{n-k}\ii \sin\pi\omega}{\pi(k-n-\omega)} \hyper{3}{2}{-n,-\omega,k-n-\omega}{-n-\omega,k-n-\omega+1}{1}=\frac{(-1)^{n-k}\ii \sin\pi\omega}{\pi(k-n-\omega)}r_{n,k}^\omega,
  \end{Eqnarray*}
  where
  \begin{displaymath}
    r_{n,k}^\omega=\hyper{3}{2}{-n,-\omega,k-n-\omega}{-n-\omega,k-n-\omega+1}{1}.
  \end{displaymath}
  It is thus sufficient to prove that $r_{n,k}^\omega=0$ for $k=0,\ldots,n-1$ and every $n\in\BB{N}$. (We already know from Corollary~1 that $r_{n,n}^\omega\neq0$, because $\SS_n^\omega$ is of degree $n$.)
  
  We do so by first proving a mixed recurrence which is of its own independent interest: the polynomials $\SS_n^\omega$, as defined by \R{eq:2.5}, satisfy
  \begin{equation}
    \label{eq:2.6}
    \SS_n^\omega(z)=z \SS_{n-1}^\omega(z)+\frac{\omega^2}{(\omega+n-1)(\omega+n)} \SS_{n-1}^{\omega-1}(z),\qquad n\in\BB{N}.
  \end{equation}
  
  At the first instance, it follows from \R{eq:2.5} that
  \begin{displaymath}
    \SS_n^\omega(z)-z\SS_{n-1}^\omega(z)=z^n\left[\sum_{\ell=0}^n (-1)^\ell \frac{(-n)_\ell (-\omega)_\ell}{\ell!(-\omega-n)_\ell} z^{-\ell} -\sum_{\ell=0}^{n-1}(-1)^\ell \frac{(-n+1)_\ell (-\omega)_\ell}{\ell!(-\omega-n+1)_\ell} z^{-\ell} \right]\!.
  \end{displaymath}
  Since
  \begin{displaymath}
    \frac{(-n)_\ell(-\omega)_\ell}{\ell!(-\omega-n)_\ell}-\frac{(-n+1)_\ell(-\omega)_\ell}{\ell!(-\omega-n+1)_\ell}=-\frac{\omega^2}{(\omega+n-1)(\omega+n)} \frac{(-n+1)_{\ell-1}(-\omega+1)_{\ell-1}}{(\ell-1)!(-\omega-n+1)_{\ell-1}},
  \end{displaymath}
  we deduce that
  \begin{Eqnarray*}
     \SS_n^\omega(z)-z\SS_{n-1}^\omega(z)&=&-\frac{\omega^2}{(\omega+n-1)(\omega+n)}z^{n-1} \!\sum_{\ell=1}^n (-1)^\ell \frac{(-n\!+\!1)_{\ell-1} (-\omega\!+\!1)_{\ell-1}}{(\ell\!-\!1)!(-\omega\!-\!n\!+\!2)_{\ell-1}} z^{-\ell+1}\\
     &=&\frac{\omega^2}{(\omega+n-1)(\omega+n)}z^{n-1} \sum_{\ell=0}^{n-1} (-1)^\ell \frac{(-n+1)_\ell (-\omega+1)_\ell}{\ell!(-\omega-n+2)_\ell} z^{-\ell}\\
     &=&\frac{\omega^2}{(\omega+n-1)(\omega+n)}z^{n-1} \hyper{2}{1}{-n+1,-\omega+1}{-\omega-n+1}{z^{-1}}\\
     &=&\frac{\omega^2}{(\omega+n-1)(\omega+n)} \SS_{n-1}^{\omega-1}(z)
  \end{Eqnarray*}
  and \R{eq:2.6} is true.
  
  Consequently, by induction, the formula being true for $n=0$ and $n=1$ by direct computation,
  \begin{Eqnarray*}
    \langle \SS_n^\omega,z^k\rangle_\omega&=& \frac{\ii}{2\pi}\int_{-\pi}^\pi \SS_n^\omega(\ee^{\ii\theta})\ee^{\ii(-k+\omega)\theta}\D\theta=\frac{\ii}{2\pi} \int_{-\pi}^\pi \SS_{n-1}^\omega(\ee^{\ii\theta})\ee^{\ii(-k+\omega)\theta}\D\theta \\
    &&\mbox{}+\frac{\omega^2}{(\omega+n-1)(\omega+n)} \frac{\ii}{2\pi} \int_{-\pi}^\pi \SS_{n-1}^{\omega-1}(\ee^{\ii\theta}) \ee^{\ii(-k+\omega-1)\theta}\D\theta\\
  &=&\langle \SS_{n-1}^\omega,z^{k}\rangle_\omega+\frac{\omega^2}{(\omega+n-1)(\omega+n)} \langle \SS_{n-1}^{\omega-1},z^{k}\rangle_{\omega-1}=0 
  \end{Eqnarray*}
  for $k=0,\ldots,n-2$.  All we need to prove is that $\langle \SS_n^\omega, z^{n-1}\rangle_\omega=0$, $\langle \SS_n^\omega,z^n\rangle\neq0$, and to this end it is sufficient to prove that $r_{n,n-1}^\omega=0$ and $r_{n,n}^\omega\neq0$ respectively. But 
  \begin{displaymath}
    r_{n,n-1}^\omega=\hyper{3}{2}{-n,-\omega,-\omega-1}{-n-\omega,-\omega}{1}=\hyper{2}{1}{-n,-\omega-1}{-n-\omega}{1}=(-1)^n \frac{(-n+1)_n}{(-n-\omega)_n}=0,
  \end{displaymath}
  where we have used the standard Vandermonde formula to sum up ${}_2F_1$ series at $z=1$. Since our stipulated form of $\SS_n^\omega$ is monic, the expression \R{eq:2.5} follows -- as does \R{eq:2.6}, which might be of an independent interest. 
  
  Finally, it follows that $r_{n,n}^\omega\neq0$ from \R{eq:2.6} by easy induction and our proof is done.
\end{proof}

\begin{corollary}
   For every $m,n\in\BB{N}$
   \begin{equation}
     \label{eq:2.7}
     z^{-n} \SS_n^m(z)=z^{-m}\SS_m^n(z).
   \end{equation}
\end{corollary}

The time has come to lift the assumption that $\omega>0$: since $\SS_n^0(z)=z^n$ is a simple and very well-known case, we need just to consider the case $\omega<0$. 

\begin{lemma} \label{lemma:neg_omega}
  Let $\omega>0$, $\omega\not\in\{1,2,\ldots,n\}$, then
  \begin{equation}
    \label{eq:2.8}
    \SS_n^{-\omega}(z)=(-1)^n \frac{(\omega)_n}{(1-\omega)_n} z^n \SS_n^{\omega-1}(z^{-1}).
  \end{equation}
\end{lemma}

\begin{proof}
  Direct substitution in \R{eq:2.5} (which has been obtained without assuming the sign of $\omega\neq0$) yields
  \begin{displaymath}
    z^n \SS_n^{-\omega}(z^{-1})=\hyper{2}{1}{-n,\omega}{-n+\omega}{-z}=\sum_{\ell=0}^n {n\choose \ell} \frac{(\omega)_{n-\ell}}{(-n+\omega)_{n-\ell}} z^{n-\ell}. 
  \end{displaymath}
  But
  \begin{Eqnarray*}
    (\omega)_{n-\ell}&=&(-1)^\ell \frac{(\omega)_n}{(-n-\omega+1)_\ell},\\
    (-n+\omega)_{n-\ell}&=&(-1)^{n-\ell} \frac{(1-\omega)_n}{(1-\omega)_\ell},
  \end{Eqnarray*}
  therefore, following simple algebra,
  \begin{displaymath}
    z^n \SS_n^{-\omega}(z^{-1})=(-1)^n \frac{(\omega)_n}{(1-\omega)_n} \SS_n^{\omega-1}(z).
  \end{displaymath}
  The expression \R{eq:2.8} follows by replacing $z$ with $z^{-1}$. 
\end{proof}

Note that $\SS_n^m$ blows up for $m\in\{1,2,\ldots,n\}$, because the denominator of the hypergeometric function vanishes. This comes as a little surprise: the main message of \R{eq:2.8} is that, flipping the sign of a non-integer $\omega$ is equivalent to conformally reflecting a skyburst polynomial (with a unit shift in parameter) in respect to the unit circle. Such a reflection takes the origin to infinity and we have already seen, e.g.\ in \R{eq:2.7}, that $\SS_n^\omega$ vanishes at the origin for $\omega\in\{0,1,\ldots,n-1\}$.

\setcounter{equation}{0}
\setcounter{figure}{0}
\section{Recurrences and generating functions}

We have already obtained the mixed recurrence relation \R{eq:2.6} in the course of proving Theorem~\ref{th:explicit_form}. This formula is interesting in the following sense. Polynomials orthogonal on the unit circle with respect to a real-valued measure obey the {\em Szeg\H{o} recurrence\/}
\begin{equation}
  \label{eq:3.1}
  p_{n+1}(z)=zp_n(z)-\bar{\alpha}_n p_n^*(z),\qquad n\in\BB{N},
\end{equation}
for a sequence of {\em Verblunski coefficients\/} $\alpha_n$ \cite[p.~2]{simon05opuc1}, where $p_n^*(z)=z^n\overline{p_n(\bar{z}^{-1})}$. It is easy, though, to compute
\begin{displaymath}
  {\SS_n^\omega}^*(z)=\hyper{2}{1}{-n,-\omega}{-n-\omega}{-z}
\end{displaymath}
and verify that \R{eq:3.1} does not hold for skyburst polynomials. This is not very surprising, since the underlying measure is complex valued. However, the surprising fact is that the above recurrence is replaced by \R{eq:2.6}: instead of a conjugate ${\SS_n^\omega}^*$ we have $\SS_n^{\omega-1}$, with a shifted parameter. Moreover, in this section we prove several other recurrence relations.

The following mixed recurrence can be proved directly.

\begin{lemma}
  For every $n\in\BB{N}$
  \begin{equation}
    \label{eq:3.2}
    \SS_n^{\omega+1}(z)=\SS_n^\omega(z)+\frac{nz^2}{(\omega+n)(\omega+n+1)} \SS_{n-1}^\omega(z).
  \end{equation}
\end{lemma}

\begin{proof}
  We compute directly, using \R{eq:2.5},
  \begin{Eqnarray*}
    \SS_n^{\omega+1}(z)-\SS_n^\omega(z)&=&\sum_{\ell=1}^{n} {n\choose\ell} \left[\frac{(-\omega-1)_\ell}{(-\omega-1-n)_\ell}-\frac{(-\omega)_\ell}{(-\omega-n)_\ell}\right] z^{n-\ell}\\
    &=&\sum_{\ell=1}^n {n\choose\ell} \frac{(-\omega)_{\ell-1}}{(-\omega-1-n)_{\ell+1}} \ell n z^{n-\ell}\\
    &=&n^2z\sum_{\ell=0}^{n-1} {{n-1}\choose\ell} \frac{(-\omega)_\ell}{(-\omega-n-1)_{\ell+1}} z^{n-1-\ell}\\
    &=&\frac{n^2z}{(\omega+n)(\omega+n+1)} \sum_{\ell=0}^{n-1}{{n-1}\choose\ell} \frac{(-\omega)_\ell}{(-\omega-n+1)_\ell} z^{n-1-\ell}\\
    &=&\frac{n^2z}{(\omega+n)(\omega+n+1)} \SS_{n-1}^\omega(z).
  \end{Eqnarray*}
\end{proof}

As a gateway to further recurrence relations, we prove the existence of a surprisingly neat generating function.

\begin{theorem}\label{th:genfunction}
  It is true that
  \begin{equation}
    \label{eq:3.3}
    \sum_{n=0}^\infty \frac{(1+\omega)_n}{n!} \SS_n^\omega(z)T^n=\frac{(1+T)^\omega}{(1-zT)^{\omega+1}},\qquad |zT|<1.
  \end{equation}
\end{theorem}

\begin{proof}
  We have
  \begin{Eqnarray*}
    G_\omega(z,T)&=&\sum_{n=0}^\infty \frac{(1+\omega)_n}{n!} \SS_n^\omega(z) T^n=\sum_{n=0}^\infty \frac{(1+\omega)_n}{n!} T^n \sum_{\ell=0}^n (-1)^\ell \frac{n!(-\omega)_\ell z^{n-\ell}}{\ell!(n-\ell)!(-\omega-n)_\ell}\\
    &=&\sum_{\ell=0}^\infty \frac{(-1)^\ell}{(-\omega)_\ell}{\ell!} \sum_{n=\ell}^\infty \frac{(1+\omega)_n z^{n-\ell}T^n}{(n-\ell)!(-\omega-n)_\ell}\\
    &=&\sum_{\ell=0}^\infty \frac{(-1)^\ell (-\omega)_\ell T^\ell}{\ell!} \sum_{n=0}^\infty \frac{ (1+\omega)_{n+\ell}(Tz)^n}{n!(-\omega-n-\ell)_\ell}.
  \end{Eqnarray*}
  But 
  \begin{displaymath}
    \frac{ (1+\omega)_{n+\ell}}{(-\omega-n-\ell)_\ell}=\frac{(1+\omega)_n (n+1+\omega)_\ell}{(-1)^\ell (\omega+n+1)_\ell}=(-1)^\ell(\omega+1)_n,
  \end{displaymath}
  therefore
  \begin{Eqnarray*}
     G(z,T)&=&\sum_{\ell=0}^\infty \frac{(-\omega)_\ell T^\ell}{\ell!} \sum_{n=0}^\infty \frac{(\omega+1)_n(Tz)^n}{n!}=\hyper{1}{0}{-\omega}{\mbox{---}}{T} \hyper{1}{0}{\omega+1}{\mbox{---}}{Tz}\\
     &=&\frac{(1+T)^\omega}{(1-Tz)^{\omega+1}},
  \end{Eqnarray*}
  summing up the binomial ${}_1F_0$ series explicitly \cite[p.~74]{rainville60sf}.
\end{proof}

The generating function \R{eq:3.3} is a pathway to a wide array of results. For example, expanding $G_\omega(-1,T)$, we obtain at once 
\begin{equation}
  \label{eq:3.4}
  \SS_n^\omega(-1)=(-1)^n \frac{n!}{(1+\omega)_n}\neq0,\qquad n\in\BB{Z}_+,\quad \omega>0,
\end{equation}
an expression which will be useful in Section~4. With minor effort, \R{eq:3.4} can be generalised.

\begin{lemma}
  \label{th:dervs}
  For every $m\in\BB{Z}_+$ it is true that
  \begin{equation}
    \label{eq:3.5}
    \frac{\D^m \SS_n^\omega(-1)}{\D z^m}=(-1)^{n-m} n! \frac{(1+\omega)_m}{(1+\omega)_n} {n\choose m},\qquad n\geq m,\quad \omega>0.
  \end{equation}
  Therefore
  \begin{displaymath}
      \SS_n^\omega(z)=\frac{(-1)^n}{(1+\omega)_n} \hyper{2}{1}{-n,1+\omega}{1}{1+z}\!.
  \end{displaymath}
\end{lemma}

\begin{proof}
  Differentiating \R{eq:3.3} and letting $z=-1$, we have 
  \begin{Eqnarray*}
    \sum_{n=m}^\infty \frac{(1+\omega)_n}{n!} \frac{\D^m \SS_n^\omega(z)}{\D z^m} T^n&=&(1+T)^\omega \frac{\D^m}{\D z^m} \frac{1}{(1-Tz)^{\omega+1}}=\frac{(1+\omega)_mT^m}{(1+T)^{m+1}}\\
    &=&(1+\omega)_m T^m \hyper{1}{0}{m+1}{\mbox{---}}{-T}\\
    &=&(1+\omega)_m \sum_{n=m}^\infty {n\choose m} (-1)^{n-m} T^n
  \end{Eqnarray*}
  and deduce \R{eq:3.5}. Therefore
  \begin{Eqnarray*}
    \SS_n^\omega(z)&=&\sum_{m=0}^n \frac{1}{m!} \frac{\D^m p_n^\omega(-1)}{\D z^m}=\frac{(-1)^n n!}{(1+\omega)_n} \sum_{m=0}^n (-1)^m {n\choose m} \frac{(1+\omega)_m}{m!}(1+z)^m\\
    &=&\frac{(-1)^nn!}{(1+\omega)_n} \hyper{2}{1}{-n,1+\omega}{1}{1+z}
  \end{Eqnarray*}
  and the proof is complete.
\end{proof}

The generating function also lends itself toward the derivation of lifting and lowering recurrences for skyburst polynomials.

\begin{theorem}\label{th:recurrences}
  The following {\em lifting\/}
  \begin{equation}
    \label{eq:3.6}
    \frac{(2+\omega)_n}{n!} \SS_n^{\omega+1}(z)=(1+z)\sum_{\ell=0}^{n-1} \frac{(1+\omega)_\ell}{\ell!} z^{n-\ell-1} \SS_\ell^\omega(z)+\frac{(1+\omega)_n}{n!} z\SS_n^\omega(z)
  \end{equation}
  and {\em lowering\/}
  \begin{equation}
    \label{eq:3.7}
    \frac{(\omega)_n}{n!} \SS_n^{\omega-1}(z)= (1+z)\sum_{\ell=0}^{n-1}(-1)^{n-\ell} \frac{(1+\omega)_\ell}{\ell!} \SS_\ell^\omega(z)+\frac{(1+\omega)_n}{n!} \SS_n^{\omega}(z)
  \end{equation}
  recurrences are valid for every $n\in\BB{Z}_+$.
\end{theorem}

\begin{proof}
  We commence by proving \R{eq:3.6}. Since
  \begin{displaymath}
      \frac{1+T}{1-Tz}=1+(1+z)\sum_{\ell=1}^\infty z^{\ell-1}T^\ell,
  \end{displaymath}
  it follows from \R{eq:3.3} that
  \begin{Eqnarray*}
    &&\sum_{n=0}^\infty \frac{(2+\omega)_n}{n!} \SS_n^{\omega+1}(z)T^n=\frac{(1+T)^{\omega+1}}{(1-Tz)^{\omega+2}}\\
    &=&\left[1+(1+z)\sum_{\ell=1}^\infty z^{\ell-1}T^\ell\right] \sum_{n=0}^\infty \frac{(1+\omega)_n}{n!} \SS_n^\omega(z)T^n\\
    &=&\sum_{n=0}^\infty \frac{(1+\omega)_n}{n!} \SS_n^\omega(z)T^n +(1+z)\sum_{\ell=1}^\infty z^{\ell-1} \sum_{n=0}^\infty \frac{(1+\omega)_n}{n!} \SS_n^\omega(z) T^{n+\ell}\\
    &=&\sum_{n=0}^\infty \frac{(1+\omega)_n}{n!} \SS_n^\omega(z)T^n +(1+z) \sum_{\ell-1}^\infty z^{\ell-1} \sum_{n=\ell}^\infty \frac{(1+\omega)_{n-\ell}}{(n-\ell)!} \SS_{n-\ell}^\omega(z)T^n\\
    &=&\sum_{n=0}^\infty \frac{(1+\omega)_n}{n!} \SS_n^\omega(z)T^n +(1+z)\sum_{n=1}^\infty \sum_{\ell=1}^n \frac{(1+\omega)_{n-\ell}}{(n-\ell)!} z^{\ell-1} \SS_{n-\ell}^\omega(z)T^n\\
    &=&\sum_{n=0}^\infty \frac{(1+\omega)_n}{n!} \SS_n^\omega(z)T^n +(1+z)\sum_{n=1}^\infty \sum_{\ell=0}^{n-1} \frac{(1+\omega)_\ell}{\ell!} z^{n-\ell-1}\SS_\ell^\omega(z) T^n
  \end{Eqnarray*}
  and the assertion \R{eq:3.6} follows by comparing powers of $T$.
  
  Similarly,
  \begin{displaymath}
    \frac{1-Tz}{1+T}=1+(1+z)\sum_{\ell=1}^\infty (-1)^\ell T^\ell,
  \end{displaymath}
  therefore
  \begin{Eqnarray*}
    &&\sum_{n=0}^\infty \frac{(\omega)_n}{n!} \SS_n^{\omega-1}(z)T^n=\frac{(1+T)^{\omega-1}}{(1-Tz)^\omega}=\frac{1-Tz}{1+T} \times \frac{(1+T)^\omega}{(1-Tz)^{\omega+1}}\\
    &=&\left[1+(1+z)\sum_{\ell=1}^\infty (-1)^\ell T^\ell\right] \sum_{n=0}^\infty \frac{(1+\omega)_n}{n!} \SS_n^\omega(z) T^n\\
    &=&\sum_{n=0}^\infty \frac{(1+\omega)_n}{n!} \SS_n^\omega(z)T^n +(1+z)\sum_{\ell=1}^\infty (-1)^\ell \sum_{n=0}^\infty \frac{(1+\omega)_n}{n!} \SS_n^\omega(z) T^{n+\ell}\\
    &=&\sum_{n=0}^\infty \frac{(1+\omega)_n}{n!} \SS_n^\omega(z)T^n +(1+z) \sum_{n=1}^\infty \sum_{\ell=0}^{n-1} (-1)^{n-\ell} \frac{(1+\omega)_\ell}{\ell!} \SS_\ell^\omega(z) T^n
  \end{Eqnarray*}
  and, comparing powers of $T$, the proof of \R{eq:3.7} follows.
\end{proof}

The generating function \R{eq:3.3} is a convenient pathway towards a differential recurrence for skyburst polynomials.

\begin{lemma}
  The recurrence
  \begin{equation}
    \label{eq:3.8}
    (\omega+n)\frac{\D\SS_n^\omega(z)}{\D z}=nz\frac{\D\SS_{n-1}^\omega(z)}{\D z}+n(1+\omega)\SS_{n-1}^\omega(z)
  \end{equation}
  holds for all $n\in\BB{Z}_+$.
\end{lemma}

\begin{proof}
  We again commence from the generating function \R{eq:3.3}. Since
  \begin{displaymath}
    \sum_{n=0}^\infty \frac{(1+\omega)_n}{n!} \frac{\D \SS_n^\omega(z)}{\D z}T^n=(\omega+1)T \frac{(1+T)^\omega}{(1-Tz)^{\omega+2}},
  \end{displaymath}
  we have
  \begin{Eqnarray*}
    (1-Tz) \sum_{n=0}^\infty \frac{(1+\omega)_n}{n!} \frac{\D \SS_n^\omega(z)}{\D z}T^n&=& (1+\omega)T \frac{(1+T)^\omega}{(1-Tz)^{\omega+1}} \\
    &=&(1+\omega)T \sum_{n=0}^\infty \frac{(1+\omega)_n}{n!} \SS_n^\omega(z)T^{n+1}.
  \end{Eqnarray*}
  Therefore, since $\SS_0^\omega\equiv1$,
  \begin{Eqnarray*}
    &&\sum_{n=1}^\infty \frac{(1+\omega)_n}{n!} \frac{\D \SS_n^\omega(z)}{\D z} T^n-z\sum_{n=1}^\infty \frac{(1+\omega)_{n-1}}{(n-1)!} \frac{\D \SS_{n-1}^\omega(z)}{\D z}T^n\\
    &=&(1+\omega)\sum_{n=1}^\infty \frac{(1+\omega)_{n-1}}{(n-1)!} \SS_{n-1}^\omega(z) T^n
  \end{Eqnarray*}
  and we deduce the differential recurrence \R{eq:3.8} comparing powers of $T$.
\end{proof}

Finally in this section, we demonstrate that the skyburst polynomials obey a second-order linear differential equations -- something that should come as little surprise, because of their relation to hypergeometric functions.

\begin{theorem}
  \label{th:ODE}
  The function $\SS_n^\omega$ obeys the differential equation
  \begin{equation}
    \label{eq:3.9}
    -z(1+z)\frac{\D^2\SS_n^\omega(z)}{\D z^2}+[1-(2+\omega-n)(z+1)]\frac{\D\SS_n^\omega(z)}{\D z}+(1+\omega)n\SS_n^\omega(z)=0,
  \end{equation}
  with regular-singular points at $-1$ and 0.
\end{theorem}

\begin{proof}
  Our starting point is Lemma~\ref{th:dervs}. A hypergeometric function
  \begin{displaymath}
    y(z)=\hyper{2}{1}{a,b}{c}{z}
  \end{displaymath}
  obeys the differential equation
  \begin{displaymath}
    x(1-x)y''(x)+[c-(a+b+1)x]y'(x)-aby(x)=0
  \end{displaymath}
  \cite[p.~54]{rainville60sf}. The equation \R{eq:3.9} follows by letting $a=-n$, $b=1+\omega$, $c=1$ and $y=1+z$.
\end{proof}

\begin{corollary}\label{cor:double_zeros}
 All the zeros of  $\SS_n^\omega$, except possibly at the origin, are simple. 
\end{corollary}

\begin{proof}
  Suppose that $\SS_n^\omega(\tilde{z})=\D\SS_n^\omega(\tilde{z})/\D z=0$ for some $\tilde{z}\in\BB{C}\setminus\{-1,0\}$. Solving \R{eq:3.9} with these initial conditions we obtain $\SS_n^\omega\equiv0$, a contradiction. Moreover, $\SS_n^\omega(-1)\neq0$ according to \R{eq:3.4}.
 \end{proof}

Note that it follows at once from \R{eq:2.5} that
\begin{displaymath}
  \SS_n^\omega(0)=\frac{n!(-\omega)_n}{(-n-\omega)_n}\neq0
\end{displaymath}
unless $\omega\in\{0,1,\ldots,n-1\}$, when \R{eq:2.7} implies that $\SS_n^\omega$ has a zero of multiplicity $n-m$.

\setcounter{equation}{0}
\setcounter{figure}{0}
\section{The pattern of the zeros}

Orthogonal polynomials on the real line with respect to complex-valued highly oscillatory measures have been investigated in \cite{celsus22kpt}. Perhaps their most fascinating feature is the behaviour of their zeros. The zeros no longer reside in the support of the measure. As an example (as a matter of fact, the only example investigated at depth), consider the inner product
\begin{displaymath}
  \langle\!\langle f,g\rangle\!\rangle_\omega=\int_{-1}^1 f(x)g(x)\ee^{\ii\omega x}\D x,\qquad \omega\geq0.
\end{displaymath}
As the parameter $\omega$ grows, the zeros of the orthogonal polynomial $p_n^\omega$ trace $n$ trajectories in the complex plane: each such trajectory commences at a zero of a Legendre polynomial and, as $\omega\rightarrow\infty$, tends to either $+1$ or $-1$. However -- and this is the reason to their name in \cite{celsus22kpt}, {\em kissing polynomials,\/} these trajectories do not stay distinct: at certain points $\omega^*$ they `kiss': the trajectory of $p_n^{\omega^*}$ briefly touches the trajectory of the zeros of $p_{n-1}^{\omega^*}$. These are precisely the points where the Hankel determinant $\GG{t}_n^{\omega^*}$ vanishes and, at that instance, $p_{n-1}^{\omega^*}$ and $p_n^{\omega^*}$ coincide. (In other words, $p_n^{\omega^*}$ is of degree $n-1$.) 

\begin{figure}[htb]
  \begin{center}
    \hspace*{-30pt}\includegraphics[width=130pt]{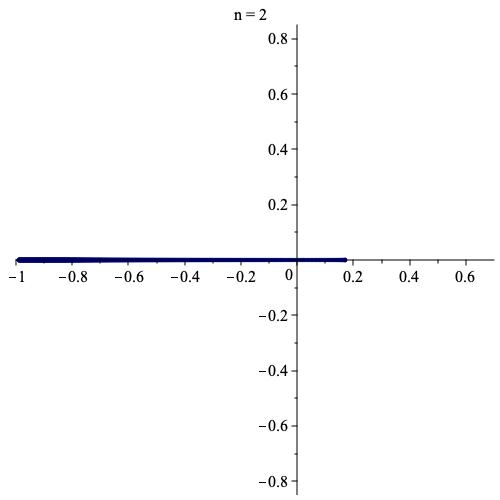}\hspace{10pt}\includegraphics[width=130pt]{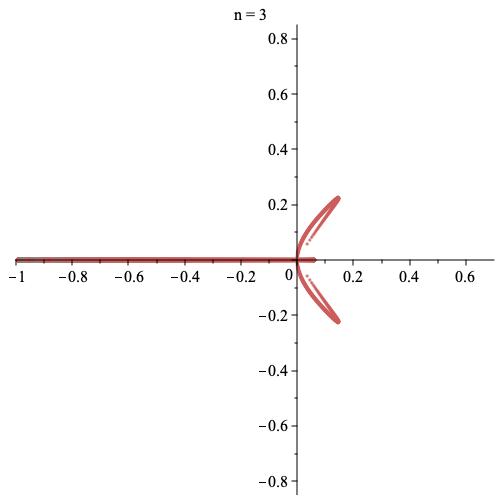}\hspace{10pt}\includegraphics[width=130pt]{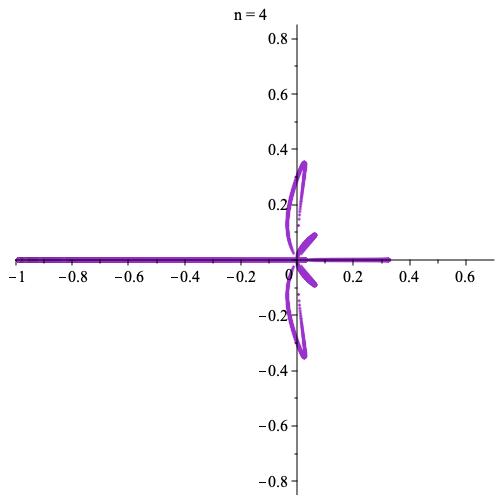}\hspace*{-30pt}
    
     \hspace*{-30pt}\includegraphics[width=130pt]{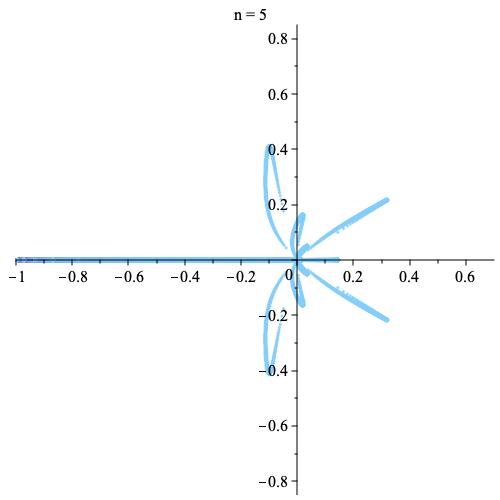}\hspace{10pt}\includegraphics[width=130pt]{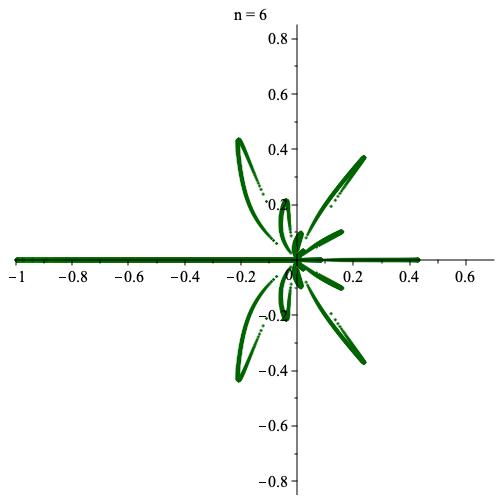}\hspace{10pt}\includegraphics[width=130pt]{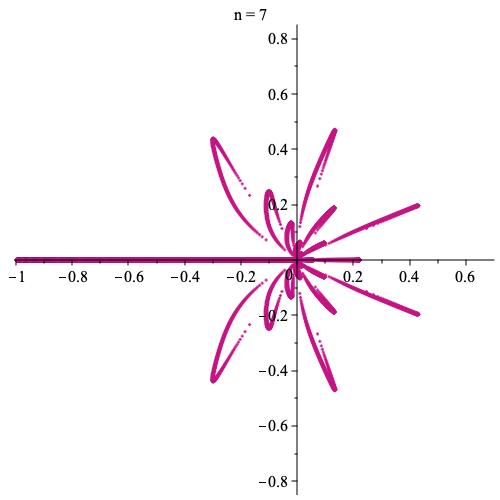}\hspace*{-30pt}
        
      \hspace*{-30pt}\includegraphics[width=130pt]{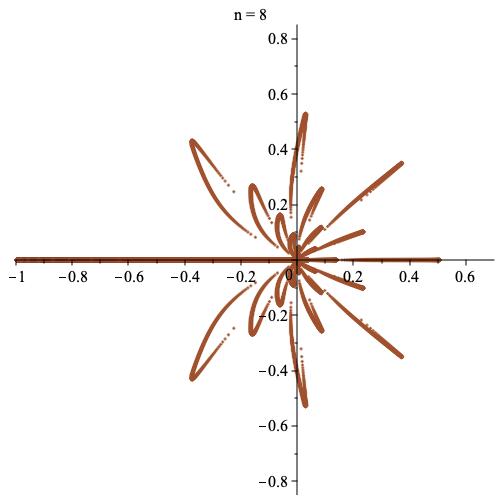}\hspace{10pt}\includegraphics[width=130pt]{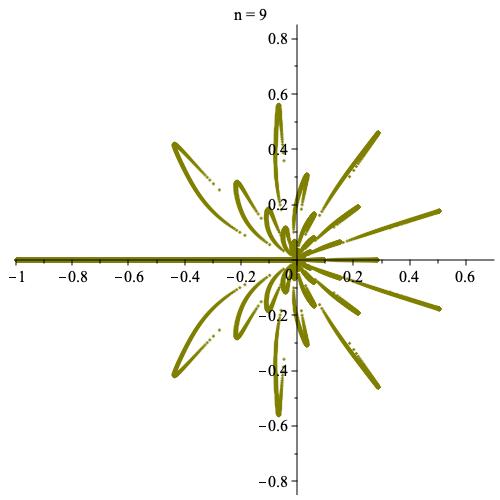}\hspace{10pt}\includegraphics[width=130pt]{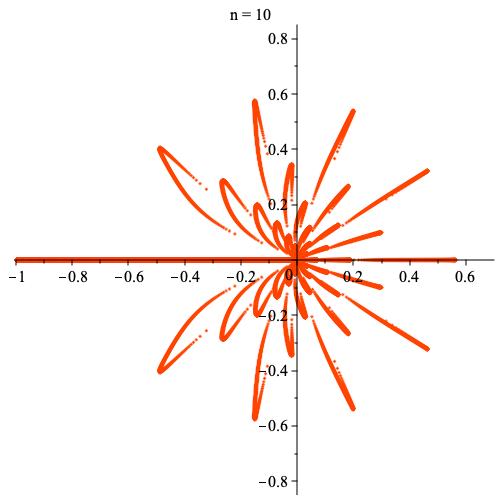}\hspace*{-30pt}
    \caption{Zeros of $\CC{S}_n^\omega$ for $n=2,\ldots,10$.}
    \label{Fig:4.1}
  \end{center}
\end{figure}

We might expect similar behaviour from the OPUC $\CC{S}_n^\omega$ defined by the bilinear form \R{eq:2.1}, yet this is not the case! We focus in this section on $\omega>0$ because, in light of Lemma~\ref{lemma:neg_omega}, the pattern for $\omega<0$ is the same, subject to the conformal map $z\rightarrow z^{-1}$ and unit shift of $\omega$.

\begin{figure}[htb]
  \begin{center}
     \includegraphics[width=260pt]{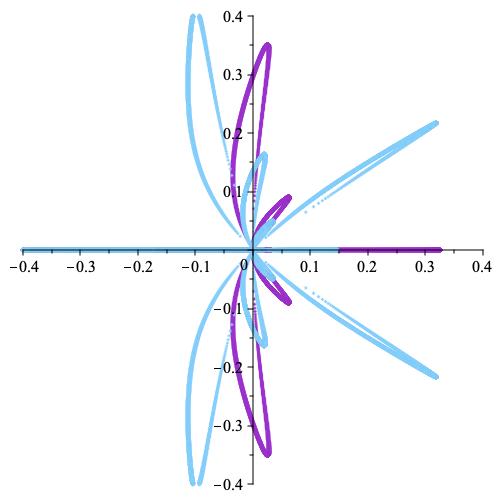}
     \caption{A closeup of the zeros of $\CC{S}_4^\omega$ and $\CC{S}_5^\omega$.}
    \label{Fig:4.2}
  \end{center}
\end{figure}

Fig.~\ref{Fig:4.1} displays the zero trajectories of $\CC{S}_n^\omega$, $2\leq n\leq 10$, in the complex plane as $\omega$ varies between 0 and $+\infty$. It provides an explanation for the name we have endowed them with, {\em skyburst polynomials.\/} As $\omega$ increases from the origin, the zeros `burst' from the origin (not surprising, since $\CC{S}_n^0(z)=z^n$) into the complex plane, after a while the $n$ trajectories all loop back to the origin, zeros become real, negative and they slowly move towards $-1$, a point they collectively reach as $\omega\rightarrow+\infty$. However, they do now kiss: the trajectories remain separate from each other. This is evident from Fig.~\ref{Fig:4.2}, a closeup of the zeros of $\CC{S}_4^\omega$ and $\CC{S}_5^\omega$. The trajectories cross each other, but these encounters occur at distinct values of $\omega$: it is possible for a zero of $\CC{S}_4^{\omega_1}$ to coincide with a zero of $\CC{S}_5^{\omega_2}$ for $\omega_1\neq\omega_2$. The reason is that, once $\GG{t}_n^\omega$ vanishes, so do other minors in  the determinantal representation of $\CC{S}_n^\omega$,
\begin{displaymath}
  \SS^\omega_n(z)=\frac{1}{\GG{t}_n^\omega} \det\!\left[
  \begin{array}{ccccc}
    \mu_0 & \mu_1 & \cdots & \mu_{n-1} & 1\\
    \mu_{-1} & \mu_0 & \cdots & \mu_n & z\\
    \vdots & \vdots & & \vdots & \vdots\\
    \mu_{-n} & \mu_{-n+1} & \cdots & \mu_{-1} & z^n
  \end{array}
  \right]\!,
\end{displaymath}
where the $\mu_n$s are the moments.

More detailed examination of the zero trajectories of $\SS_n^\omega$ for $n\geq2$ reveals an intriguing pattern. We claim -- and this will be proved rigorously in the sequel -- that, as $\omega>0$ grows, there are two regimes:
\begin{description}
\item[ Burst:]  For $m=1,\ldots,n-1$ and $m-1\leq \omega\leq m$ the pattern is as follows: at $\omega=m-1$ $n-m$ trajectories `burst' into the complex plane at angles which are multiples of $2\pi/(n-m)$: if $n-m$ is odd, one of them does so along the positive ray. As $\omega$ grows, the trajectories sketch a loop, ultimately returning to the origin when $\omega=m$. The trajectory along the positive ray (if it exists) is `flattened': at certain point in $(m-1,m)$ it loops back to the origin, all along positive values. The remaining $m$ zeros of $\SS_n^\omega$ live in $(-1,0)$.
\item[Fizzle:] The fireworks are over once $\omega>n$. All the zeros are then `trapped' in $(-1,0)$, ultimately tending to $-1$ as $\omega\rightarrow\infty$. 
\end{description}

We now confirm the above claims. Since $\SS_1^\omega(z)=z+\frac{\omega}{1+\omega}$, for $n=1$ we have a single zero in $(-1,0)$. Let us now progress by induction, assuming that all the zeros of $\SS_m^n$ are in $(-1,0)$ for $m\geq n+1$. Because of \R{eq:2.7}, it thus follows that $\SS_n^m(z)=z^{n-m}\SS_m^n(z)$, $m\in\{0,1,\ldots,n-1\}$, has zeros of multiplicity $m$ zeros -- actually, because of Corollary~3, $m$ simple zeros --  in $(-1,0)$, as well as a zero at the origin of multiplicity $n-m$. Consider $\omega=m+\varepsilon$ for $0<\varepsilon\ll1$. As $\varepsilon$ grows away from zero, $\SS_n^{m+\varepsilon}(z)=\SS_n^m(z)+\varepsilon \partial \SS_n^m(z)/\partial\varepsilon+\O{\varepsilon^2}$. But it follows from \R{eq:2.5} and \R{eq:2.7} that
\begin{displaymath}
  \SS_n^{m+\varepsilon}(z)=\frac{n!^2 z^{n-m}}{(n-m)!(n+m)!} [1+\O{z}]-\varepsilon (-1)^{n-m} \frac{m!^2(n-m-1)!}{(n+m)!}[1+\O{z}]+\O{\varepsilon^2}.
\end{displaymath}
We are interested in the $n-m$ zeros emanating from the origin for $m+\varepsilon$: suppose that such a zero is $r_\varepsilon \ee^{\ii\theta_\varepsilon}$, where $0<r\ll1$ -- actually, $r\approx \varepsilon^{1/(n-m)}$. Therefore
\begin{displaymath}
  \theta_\varepsilon\in\left\{\pi+\frac{2\pi k}{n-m}\,:\, k=0,\ldots,n-m-1\right\}
\end{displaymath}
-- the zeros leave the origin in $n-m$ trajectories. One (for $k=0$) leads into $(-1,0)$, the other emerge into $\BB{C}\setminus(-1,0)$ in $n-m-1$ equiangular rays. (If $n-m$ is even one of these rays proceeds along the positive ray.) 

Let us consider what happens as $\varepsilon$ increases in $(0,1)$. The zero in $(-1,0)$ cannot return to the origin or cross $-1$, because it follows at once from \R{eq:2.5} that $\SS_n^\omega(0)\neq0$ for non-integer $\omega$, while $\SS_n^\omega(-1)\neq0$ because of \R{eq:3.4}. There is another theoretical possibility for zeros to leave $(-1,0)$: the trajectories of two (or more) zeros coalesce at a point and thence emerge into $\BB{C}\setminus(-1,0)$. This is rules out by Corollary~\ref{cor:double_zeros} because a point of zero coalescence is a zero of nontrivial multiplicity. 

What about the $n-m-1$ zero trajectories that have emerged into $\BB{C}\setminus(-1,0)$? They evolve there as $\omega$ increases until $\omega=m+1$, when they must return to the origin because $\CC{S}_n^{m+1}(z)=z^{n-m-1}\SS_{m+1}^n(z)$ and each forms a loop because of the continuity of (simple!) zeros. No zero trajectory may approach the negative ray (in particular the interval $(-1,0)$) because the trajectory have $z\leftrightarrow\bar{z}$ symmetry and such an encounter will bring two trajectories together: again, this is ruled out by Corollary~\ref{cor:double_zeros}. Once $n-m$ is even, one of the trajectories at $\omega=m$ emerges into the positive ray and (again, by simplicity of zeros) it stays there until $\omega=m+1$. 

In other words, in each interval $\omega\in(m,m+1)$ for $m=0,1,\ldots,n-1$ there are $n-m-1$ zero trajectories looping in $\BB{C}\setminus(-1,0)$, while the remaining $m+1$ zeros live in $(-1,0)$. Once $\omega$ exceeds $n$, all the zeros are in $(-1,0)$ and it is trivial to note that, by \R{eq:2.5},
\begin{displaymath}
  \lim_{\omega\rightarrow\infty} \SS_n^\omega(z)=(1+z)^n
\end{displaymath}
-- ultimately, all the zeros congregate at $-1$. This completely explains the pattern visible in Fig.~\ref{Fig:4.1}.

\section*{Acknowledgements}

The work of the first author is supported by the project PID2021-124472NB-I00 funded by MCIN/AEI/ 10.13039/501100011033 and by ``ERDF A way of making Europe'' and the project E48\_23R from Diputaci\'on General de Arag\'on (Spain) and ERDF ``Construyendo Europa desde Arag\'on''.

\bibliographystyle{agsm}
\bibliography{SkyBurst}

@article {cantero16fop,
    AUTHOR = {Cantero, Mar\'{\i}a Jos\'{e} and Iserles, Arieh},
     TITLE = {From orthogonal polynomials on the unit circle to functional
              equations via generating functions},
   JOURNAL = {Trans. Amer. Math. Soc.},
  FJOURNAL = {Transactions of the American Mathematical Society},
    VOLUME = {368},
      YEAR = {2016},
    NUMBER = {6},
     PAGES = {4027--4063},
      ISSN = {0002-9947},
   MRCLASS = {42C05 (33C05 39A05)},
  MRNUMBER = {3453364},
MRREVIEWER = {Christian Remling},
       DOI = {10.1090/tran/6454},
       URL = {https://doi.org/10.1090/tran/6454},
}

@article {celsus22kpt,
    AUTHOR = {Celsus, Andrew F. and Dea\~{n}o, Alfredo and Huybrechs, Daan and
              Iserles, Arieh},
     TITLE = {The kissing polynomials and their {H}ankel determinants},
   JOURNAL = {Trans. Math. Appl.},
  FJOURNAL = {Transactions of Mathematics and its Applications. A Journal of
              the IMA},
    VOLUME = {6},
      YEAR = {2022},
    NUMBER = {1},
     PAGES = {tnab005, 66},
   MRCLASS = {33C47 (30C50)},
  MRNUMBER = {4366067},
MRREVIEWER = {\'{A}rp\'{a}d Baricz},
       DOI = {10.1093/imatrm/tnab005},
       URL = {https://doi.org/10.1093/imatrm/tnab005},
}

@Book{rainville60sf,
  author = 	 {E. D. Rainville},
  title = 	 {Special Functions},
  publisher = 	 {Macmillan},
  year = 	 {1960},
  address = 	 {New York}
 }

@book {simon05opuc1,
    AUTHOR = {Simon, Barry},
     TITLE = {Orthogonal polynomials on the unit circle. {P}art 1},
    SERIES = {American Mathematical Society Colloquium Publications},
    VOLUME = {54},
      NOTE = {Classical theory},
 PUBLISHER = {American Mathematical Society, Providence, RI},
      YEAR = {2005},
     PAGES = {xxvi+466},
      ISBN = {0-8218-3446-0},
   MRCLASS = {42-02 (30C85 33C45 42C05 47B36 47N50)},
  MRNUMBER = {2105088},
MRREVIEWER = {P. L. Duren},
       DOI = {10.1090/coll054.1},
       URL = {https://doi.org/10.1090/coll054.1},
}

@article {Schechter59oic,
    AUTHOR = {Schechter, Samuel},
     TITLE = {On the inversion of certain matrices},
   JOURNAL = {Math. Tables Aids Comput.},
  FJOURNAL = {Mathematical Tables and other Aids to Computation},
    VOLUME = {13},
      YEAR = {1959},
     PAGES = {73--77},
      ISSN = {0891-6837},
   MRCLASS = {65.00},
  MRNUMBER = {105798},
MRREVIEWER = {I. R. Savage},
}

\end{document}